
\magnification=1200
\input amssym.def
\overfullrule0pt
\parindent=22pt

\def\R{{\Bbb R}}

\def\C{{\Bbb C}}

\def\T{{\Bbb T}}

\def\11{{\bf 1\!\!}}
\def\ratop{\mathop{\hbox to .25in{\rightarrowfill}}\limits}

\def\nhang{\hangindent=4pc\hangafter=1}

\centerline{\bf Extension of holomorphic bundles to the disc}
\medskip

\centerline{\bf (and Serre's Problem on Stein bundles).}

\centerline{
by  Jean-Pierre Rosay%
\footnote*{Partly supported by NSF grant.}}

\footnote{} {AMS classification: 32L05, 32Q28.}

\bigskip
 \hglue3truein {\it In memory of Alain Dufresnoy, a friend

 \hglue2.5truein and an inspiring and generous Mathematician.}
\bigskip
\noindent {\tenbf Introduction.} We shall study (locally trivial) holomorphic bundles
over open sets in $\C$ (in which, the variable will be denoted by $\zeta$), 
with fiber $\C^n$. We shall assume that, in appropriate local trivializations,
the bundle is given by gluing (transition)
fiber automorphisms that are locally independent of
the base point and that belong to a group $\bf G$ of automorphisms of $\C^n$.
\smallskip\noindent
Of course, such a bundle is trivial over any simply connected region. If $U$ is
an open set in $\C$ covered by simply connected domains $(U_j)$, over $U$ any such 
fiber bundle is
obtained by gluing the trivial bundles $U_j\times \C^n$, by using fiber automorphisms
belonging to $\bf G$ and locally independent of the base point.
\smallskip\noindent 
Our goal is to extend such bundles to bundles defined over larger open sets in $\C$,
and in fact over the whole Riemann sphere. For the extended bundles,
we must allow the fiber automorphisms to (even locally) depend on the base point.
\bigskip
\bigskip
\noindent {\tenbf Definition.} We shall say that a group $\bf G$ of automorphisms of
$\C^n$ is {\it generated by one parameter groups} if there are one (complex)
parameter groups $(S^t)_{t\in \C}$ (i.e. $S^{s+t}=S^s\circ S^t$) of automorphisms of
$\C^n$ whose elements generate $\bf G$.
\bigskip\noindent
Examples to have in mind are:\bigskip\noindent
(A) The group of all polynomial automorphisms of $\C^2$.
\bigskip\noindent
(B) For any $n$, the group of polynomial automorphisms generated by the invertible affine
maps and the polynomial shears, i.e. maps of the type
$$(z_1,\cdots ,z_n)\mapsto
(z_1,\cdots ,z_{n-1},z_n+Q(z_1,\cdots ,z_{n-1}))~,$$
where $Q$ is an arbitrary polynomial in $(n-1)$
variables. A corresponding one parameter group is
$$(z_1,\cdots ,z_n)\mapsto S^t(z_1,\cdots ,z_n)=
(z_1,\cdots ,z_{n-1},z_n+t Q(z_1,\cdots ,z_{n-1}))~.$$
If $n=2$, this group (generated by affine maps and shears) is the group of
all polynomial automorphisms of $\C^2$ , by Jung's Theorem [11].
\bigskip\noindent
(C) The group of automorphisms of $\C^n$ generated by the invertible affine maps,
the shears (as in (B) but with $Q$ entire) and the `over-shears' (in the terminology
of Anders\'en-Lempert), i.e. maps of the type
$$(z_1, \cdots ,z_n)\mapsto 
(z_1,\cdots z_{n-1},z_ne^{Q (z_1,\cdots ,z_{n-1})})~.$$
In that case, set
$$S^t(z_1,\cdots ,z_n)=
(z_1, \cdots z_{n-1},z_ne^{tQ (z_1,\cdots ,z_{n-1})})~.$$
\bigskip\bigskip
\noindent One parameter groups of automorphisms of $\C^n$ are given by the flow
of complete (i.e. integrable for all time) holomorphic vector fields on $\C^n$.
There are several notions of completeness for holomorphic vector fields: in positive time, 
in real time or in complex time. In $\C^n$, these three notions coincide
[1]. On the topic of Example (C) and of complete holomorphic vector fields, see
[3] and [9].
\bigskip\bigskip\noindent
{\tenbf Proposition.}
{\it Let $D$ be a bounded domain in $\C$ bounded by a Jordan curve, and for $j=1,\cdots N$, let
$D_j$ be a region bounded by a Jordan curve. Assume that $\overline D_j \subset D$,
$\overline D_j\cap\overline D_k = \emptyset$, if $j\neq k$. Let $\Lambda = D\setminus
\cup_{j=1}^N\overline D_j$ . Let $S$  denote the Riemann sphere $S=\C\cup \{ \infty\}$.
\smallskip\noindent
Let $\Pi :~X\to \Lambda$ be a (locally trivial) holomorphic fiber bundle with
fiber $\C^n$ and (in appropriate local trivializations)
with gluing automorphisms of the fibers
locally independent of the base point and in a group $\bf G$ of automorphisms
of $\C^n$ generated by one parameter groups.
Then there exists a holomorphic fiber bundle $\tilde X \to S$ with fiber $\C^n$ and
(in appropriate local trivializations)
with gluing automorphisms of the fibers in $\bf G$ (allowed to depend on the base point),
whose restriction over $\Lambda$ is $X$.}
\bigskip\bigskip
\noindent {\bf Corollary.} {\it \hfill\break
(a) The fiber bundle in Skoda's original example of non-Stein bundle [12] 
(and see [13]) extends to a holomorphic bundle over
the Riemann sphere (and so to a holomorphic bundle over the disc).
\smallskip\noindent
(b) The possibly non-Stein fiber bundles constructed by Demailly [6] on annuli,
with fiber $\C^2$ and the gluing automorphism 
$$(z_1,z_2)\mapsto (z_2,-z_1+z_2^k)~,$$
extend to holomorphic bundles on the Riemann sphere with polynomial automorphisms
of the fibers.}
\bigskip\noindent
The Corollary is immediate, and the consequences for the problem of Stein bundles
are as follows. \hfill\break
($\alpha$) There are two versions of Skoda's example. In the first, more transparent, one
(Theorem in [12]), the base is made of the unit disc with 8 discs removed.
In the second one (final remark in [12] and Theorem 1 in [13]) the base is the disc
with 2 discs removed. In the first example, around each deleted disc the bundle is
defined by a gluing automorphism which, up to a possible permutation of the coordinates,
is an over-shear as in Example C. In the second example, around one of the deleted disc
the gluing is given by the automorphism $(z_1,z_2)\mapsto (z_1,z_2e^{z_1})$, and for the
other deleted disc by $(z_1,z_2)\mapsto (iz_2,z_1)$.
\smallskip\noindent
So,
starting from Skoda's example, one gets
a new construction of a non-trivial, non-Stein,
holomorphic bundle on the unit disc (and on $\C$). This construction 
seems to me to be much less technical than the original 
example of Demailly [6] and [7]. Furthermore the construction relies only on the
simplest case in the proof of Section II, with a trivial filling of each hole,
since around each hole the gluing is made by automorphisms that are elements of
a one parameter group of automorphisms
\smallskip\noindent
($\beta$) The maps used by Demailly are of Henon type. The map
$(z_1,z_2)\mapsto (z_2,-z_1+z_2^k)$ is the composition of
the linear map $(z_1,z_2)\mapsto (z_2,-z_1)$, and of the shear
$(z_1,z_2)\mapsto (z_1,z_2+z_1^k)$. Each of these two maps is an element of
a one parameter group of polynomial automorphisms. But, if $k\geq 2$, the map 
$(z_1,z_2)\mapsto (z_2,-z_1+z_2^k)$ itself is not even an element of a one
parameter group of (possibly non-polynomial) automorphisms of $\C^2$, since 
the proof of Case 1 in Section II
would then show that the Demailly bundles would be trivial. Following Section II,
the extension of the bundle from the annulus to the disc is then obtained by first
extending the bundle to the disc with 2 (smaller) holes, and then extending trivially
to each hole.
\smallskip\noindent
Starting from the polynomial examples of Demailly (non-Stein for a given annulus
if $k$ is large enough, and non-Stein for a given $k\geq 2$ if the annulus is thick enough),
one gets examples of  holomorphic bundles on the unit disc (and on $\C$), 
with fiber $\C^2$, and with gluing polynomial
automorphisms, that are non trivial and non-Stein. This answers a question in [6]. The question
was asked again recently by H. Skoda, who considered it to be the last question in the Serre
Problem on Stein bundles, left open after the counterexamples of Skoda [12] [13], 
Demailly [6] [7], and Coeur\'e-Loeb [5].\bigskip\noindent
If one takes $k=2$, all the polynomial automorphisms to be used in the proof of
the Proposisiton can be chosen of degree $\leq 2$.
So, there are examples of non-Stein bundles on the disc, with polynomial gluing
automorphisms of the fiber of degree $\leq 2$, while all bundles on the disc with
affine gluing automorphisms are trivial.
\vskip.5truein
\noindent {\bf II. Proof of the Proposition.}
\hfill\break
It is of course enough to extend the given bundle $X$ over each ``hole'' $D_j$,
and over the component $D_0$ of $S\setminus \overline D$, containing
${\infty}$. Fix $j\in
\{0,1,\cdots ,N\}$, and let $V$ be a neighborhood of $bD_j$ in $\Lambda$ that is
topologically an annulus. We shall decompose $V$ into the union $V^+\cup V^-$ of two
simply connected regions with intersection $V^+\cap V^-=\omega\cup\omega'$,
where $\omega$ and $\omega'$ are disjoint connected open sets. We wish to extend the bundle
$X|\Pi^{-1}(V)$ to a bundle over $\overline D_j\cup V$, defined by automorphisms of
the fibers in the group $\bf G$. Un-needed for the easy Case 1 below, un-necessary
but helpful
for the general case, we shall be more specific. For $r>0$, let $\Delta_r$ be the open disc in
$\R^2$ of radius $r$ centered at the origin.  We can take $V$ and $V^\pm$ so that 
for some homeomorphism $\chi$ of  
$\Delta_2$ into S: $D_j=\chi (\Delta_1)$, and $V=\chi (\Delta_2\setminus
\overline \Delta_1)$, $V^+$ and $V^-$ being respectively the image of
the regions ${\rm Im}~\zeta >-{1\over 2}$ and
${\rm Im}~\zeta <{1\over 2}$.\medskip\noindent
Since we start with locally constant gluing fiber automorphisms, and since 
$V^+$ and $V^-$ are simply connected, the fiber bundle $X$ is trivial over $V^+$
and over $V^-$. Over $V$, it can be defined by gluing $V^+\times \C^n$ and
$V^-\times \C^n$ over $\omega \cup \omega'$, with the identification of
$(\zeta ,z )\in V^+\times \C^n$  \hfill\break
with $(\zeta ,T_0z)\in V^-\times \C^n$ if $\zeta\in \omega$, and
\hfill\break
with $(\zeta ,T_1z)\in V^-\times \C^n$ if $\zeta\in \omega'$, \hfill\break
where $T_0$ and $T_1$ are automorphisms of $\C^n$ that belong to $\bf G$. 
Of course we could take $\T_0={\bf 1}$ (but we don't, in order to
simplify the exposition in the second case below).
We now split the proof in two cases.
\bigskip\noindent
1) Case 1. This is the easy case when there is in $\bf G$ a one (complex) parameter
group of automorphisms $(S^t)_{t\in \C}$ such that $T_1\circ T_0^{-1}=S^1$.
In that case, the bundle extends trivially over $D_j$, and the trivialization is by
means of fiber automorphisms in $\bf G$. (If one wished to get a version of
the Proposition for bundles defined by gluing automorphisms depending
on the base point, one would need here that both $T_0$ and
$T_1$ be defined
on $V^+\cup V^-$, as holomorphic function of the base point, 
and the possibity of a holomorphic choice of the
one parameter group).

\bigskip\noindent
For identifying the trivial bundle $(\overline D_j\cup V)\times \C^n$
with the restriction of $X$ over $V$, we need to define fiber automorphisms:
$$\Phi^+: V^+\times \C^n \to V^+\times \C^n$$
and
$$\Phi^-: V^-\times \C^n \to V^-\times \C^n~~,$$
such that
$$(\ast )~~~(\Phi^-)^{-1}\circ \tilde T_0\circ \Phi^+={\bf 1}~~{\rm on} ~\omega\times\C^n~,$$
$$(\ast\ast)~~~(\Phi^-)^{-1}\circ \tilde T_1\circ \Phi^+={\bf 1}~~{\rm on} ~\omega'\times\C^n~,$$
where $\tilde T_j(\zeta ,z)=(\zeta ,T_j(z))$.
\hfill\break
Let $L^+$ and $L^-$ be holomorphic functions defined respectively on
$V^+$ and $V^-$, such that\hfill\break
$L^++L^-=0$ on $\omega$\hfill\break
$L^++L^-=-1$ on $\omega'$.\hfill\break
This is a classical Cousin problem that here is elementary to solve
(use logarithms!).
Set
$$\Phi^+(\zeta , z)=\big(\zeta ,T_0^{-1}S^{L^+(\zeta )}(z)\big)$$
$$\Phi^-(\zeta ,z)=\big(\zeta , S^{-L^-(\zeta )}(z)\big)~.$$
(Here comes the local dependence on $\zeta$.)
\medskip\noindent
For $\zeta \in \omega$:
$$[S^{-L^-(\zeta )}]^{-1}\circ T_0\circ (\T_0^{-1}\circ S^{L^+(\zeta )})
= S^{L^-(\zeta )+ L^+(\zeta )} = S^0={\bf 1}~,$$
so $(\ast )$ is satisfied.\hfill\break
For $\zeta \in \omega'$:
$$[S^{-L^-(\zeta )}]^{-1}\circ T_1\circ (\T_0^{-1}\circ S^{L^+(\zeta )})
= S^{L^-(\zeta )}\circ S^1 \circ S^{L^+(\zeta )} = 
S^{L^-(\zeta )+1+L^+(\zeta )}=S^0={\bf 1}~,$$
so $(\ast\ast )$ is satisfied.\hfill\break
That ends the proof of Case 1. Notice that in case a trivial extension is
possible over each $D_j$, it does not give a global trivial extension, even just on
$\C$ or $D$.
Skoda's example illustrates that.
\bigskip\bigskip
\noindent 2) General Case. The general case can be reduced to the previous case
by introducing a higher connectivity (replacing $D_j$ by finitely many
smaller holes) before filling each hole individually
trivially. Recall that the `hole' $D_j$ is parameterized by $\Delta_1$, under the map
$\chi$, and that $\overline D_j\cup V$ is parameterized by $\Delta_2$. In  general one has
$$T_1\circ T_0^{-1}=E_k\circ \cdots \circ E_1~,$$ where each $E_p$ 
($p=1,\cdots ,k$) is an element of
a one parameter group of automorphisms in $\bf G$.
\medskip\noindent
Fix points
$$b_0=-2~,~a_1=-1<b_1<a_2<b_2<\cdots <a_p<b_p<a_{p+1}<\cdots<b_k=1~,~a_{k+1}=2~.$$
Let $W_0$ be the complement in $\Delta_2$ of the circles of
diameter $[a_p,b_p]$, $p=1,\cdots ,k$, and let $W = \chi (W_0)$.
So, $W$ is a multiply connected region in S that contains
$V$ and that is contained in $V\cup \overline D_j$. 
Set $W^+=\chi (W_0\cap \{{\rm Im}~\zeta \geq 0\})$, and
    $W^-=\chi (W_0\cap \{{\rm Im}~\zeta \leq 0\})$.
We first extend the restriction of the given fiber bundle $X$ over $V$,
to a fiber bundle $X'$ over $W$. For that purpose, consider the fiber bundle 
$X'$ obtained by gluing $W^+\times \C^n$ with $W^-\times \C^n$ over each
component $\chi ([b_p,a_{p+1}])$ of the intersection, by using the
automorphism $E_p\circ \cdots \circ E_1\circ T_0$ ($T_0$ if $p=0$, 
$T_1$ if $p=k$). Now, for the extended bundle $X'$, over each hole of $W$ 
(the images under $\chi$ of the discs of diameter $[a_p,b_p]$), we are in Case 1,
since 
$$[E_p\circ \cdots E_1\circ T_0] \circ
[E_{p-1}\circ \cdots E_1 \circ T_0]^{-1}=E_p~.$$ 
For $p=1$, read $[E_1\circ T_0]\circ T_0^{-1}=E_1$.
Therefore, according to Case 1, $X'$ extends to a bundle over $\Lambda \cup
\overline D_j$, with gluing automorphisms in $\bf G$ (trivially over
each hole of $W$, but not trivially over $D_j$ as illustrated by Demailly's
example).
\vskip1truein\noindent
{\bf Remarks.} It has been pointed out that the map 
$(z_1,z_2)\mapsto (z_2,-z_1+z_2^k)$, for $k\geq 2$, is not an element of
a one parameter group of automorphisms of $\C^2$, i.e. it is not the time-1
map of a complete vector field on $\C^2$. The reason given was rather indirect.
For a more direct approach to this question, see [2] and [4]. F. Forstneri\v c pointed out
to me the interesting fact that, although they are far from being trivial since they
can even be non-Stein, all holomorphic fibers bundles with fiber $C^n$ over a Stein base have
holomorphic sections [10], [8], and they satisfy the ``Oka Principle''.
\vskip.5truein

\centerline{\bf REFERENCES}
\bigskip
\item{[1]} P.\ Ahern, M.\ Flores, J-P.\ Rosay.
On $\R^+$ and $\C$ complete holomorphic vector fields, Proc. A.M.S. 128 (2000),
3107-3113.
\smallskip

\item{[2]} P.\ Ahern, F.\ Forstneri\v c. One parameter automorphism groups on $\C^2$,
Complex Varables Theory Appl. 27 (1995), 245-268.
\smallskip

\item{[3]} E.\ Anders\'en, L.\ Lempert. On the group of holomorphic
automorphisms of $\C^n$, Inven. Math. 110 (1992), 371-388.
\smallskip

\item{[4]} G.\ Buzzard, J.E.\ Fornaess. Complete holomorphic vector fields and time-1 maps.
Indiana Math. J. 44 (1995), 1175-1182.

\smallskip

\item{[5]} G.\ Coeur\'e, J-J.\ Loeb. A counterexample to the Serre problem with
a bounded domain in $\C^2$, Ann. Math. 122 (1985), 329-334.
\smallskip

\item{[6]} J-P.\ Demailly. Diff\'erents exemples de fibr\'es holomorphes
non de Stein. {\it Seminaire P.Lelong, H. Skoda 1976/77}, Springer L.N. in Math.
694 (1978), 15-41. 
\smallskip

\item{[7]} J-P.\ Demailly. Un exemple de fibr\'e holomorphe non de Stein
\`a fibre $\C^2$ ayant pour base le disque ou le plan, Inven. Math. 48 (1978),
293-302.
\smallskip

\item{[8]} F.\ Forstneri\v c, J.\ Prezelj. Oka's principle for holomorphic bundles with sprays,
Math. Ann. 317 (2000), 117-154.
\smallskip 

\item {[9]} F.\ Forstneri\v c, J-P.\ Rosay. Approximation of biholomorphic mappings
by automorphisms of $\C^n$, Inven. Math. 112 (1993), 323-349. Erratum in
Inven. Math. 118 (1994), 573-574.
\smallskip

\item{[10]} M.\ Gromov, Oka's principle for holomorphic sections of elliptic bundles,
Journal A.M.S. 2 (1989), 851-897.
\smallskip

\item{[11]} H. Jung. \"Uber ganze birationale Transformationen der Ebene,
Reine Angew. Math. 184 (1942), 161-174.

\smallskip

\item{[12]} H.\ Skoda. Fibr\'es holomorphes \`a base et fibre de Stein,
C.R.Acad.Sci. Paris Serie AB 284 (1977) $n^o$19 A1199-A1202.
\smallskip

\item{[13]} H.\ Skoda. Fibr\'es holomorphes \`a base et \`a fibre de Stein,
Inven. Math. 43 (1977), 97-107. 
\smallskip
\bigskip
\centerline{- - - - - - - - - - - - - - - - - - - - - }
\bigskip

\nhang{J-P.\ Rosay: Department of Mathematics, University of Wisconsin,
Madison WI 53706 USA. {\it jrosay@math.wisc.edu}}

\bye